\documentclass[11pt]{amsart}

\pdfoutput=1
\usepackage[letterpaper,margin=1in,footskip=0.25in]{geometry}
\usepackage[alphabetic,initials]{amsrefs}
\usepackage{amsmath,amsfonts,amsthm,mathrsfs}
\usepackage{amssymb}
\usepackage{mathtools}

\usepackage{hyperref}
\usepackage[usenames,dvipsnames]{xcolor}
\hypersetup{colorlinks=true,citecolor=NavyBlue,linkcolor=BrickRed,urlcolor=Orange}

\usepackage{enumitem}

\usepackage{chngcntr}

\usepackage{tikz}

\usepackage{tikz-cd}

\usetikzlibrary{matrix,arrows}
\newlength{\myarrowsize} 

\pgfarrowsdeclare{cmto}{cmto}{
	\pgfsetdash{}{0pt} 
	\pgfsetbeveljoin 
	\pgfsetroundcap 
	\setlength{\myarrowsize}{0.6pt}
	\addtolength{\myarrowsize}{.5\pgflinewidth}
	\pgfarrowsleftextend{-4\myarrowsize-.5\pgflinewidth} 
	\pgfarrowsrightextend{.8\pgflinewidth}
}{
	\setlength{\myarrowsize}{0.6pt} 
  	\addtolength{\myarrowsize}{.5\pgflinewidth}  
	\pgfsetlinewidth{0.5\pgflinewidth}
	\pgfsetroundjoin
	\pgfpathmoveto{\pgfpoint{1.5\pgflinewidth}{0}}
	\pgfpatharc{-109}{-170}{4\myarrowsize}
	\pgfpatharc{10}{189}{0.58\pgflinewidth and 0.2\pgflinewidth}
	\pgfpatharc{-170}{-115}{4\myarrowsize+\pgflinewidth}
	\pgfpathclose
	\pgfusepathqfillstroke
	\pgfpathmoveto{\pgfpoint{1.5\pgflinewidth}{0}}
	\pgfpatharc{109}{170}{4\myarrowsize}
	\pgfpatharc{-10}{-189}{0.58\pgflinewidth and 0.2\pgflinewidth}
	\pgfpatharc{170}{115}{4\myarrowsize+\pgflinewidth}
	\pgfpathclose
	\pgfusepathqfillstroke
	\pgfsetlinewidth{2\pgflinewidth}
}

\pgfarrowsdeclare{cmonto}{cmonto}{
	\pgfsetdash{}{0pt} 
	\pgfsetbeveljoin 
	\pgfsetroundcap 
	\setlength{\myarrowsize}{0.6pt}
	\addtolength{\myarrowsize}{.5\pgflinewidth}
	\pgfarrowsleftextend{-4\myarrowsize-.5\pgflinewidth} 
	\pgfarrowsrightextend{.8\pgflinewidth}
}{
	\setlength{\myarrowsize}{0.6pt} 
  	\addtolength{\myarrowsize}{.5\pgflinewidth}  
	\pgfsetlinewidth{0.5\pgflinewidth}
	\pgfsetroundjoin
	\pgfpathmoveto{\pgfpoint{1.5\pgflinewidth}{0}}
	\pgfpatharc{-109}{-170}{4\myarrowsize}
	\pgfpatharc{10}{189}{0.58\pgflinewidth and 0.2\pgflinewidth}
	\pgfpatharc{-170}{-115}{4\myarrowsize+\pgflinewidth}
	\pgfpathclose
	\pgfusepathqfillstroke
	\pgfpathmoveto{\pgfpoint{1.5\pgflinewidth}{0}}
	\pgfpatharc{109}{170}{4\myarrowsize}
	\pgfpatharc{-10}{-189}{0.58\pgflinewidth and 0.2\pgflinewidth}
	\pgfpatharc{170}{115}{4\myarrowsize+\pgflinewidth}
	\pgfpathclose
	\pgfusepathqfillstroke
	\pgfpathmoveto{\pgfpoint{1.5\pgflinewidth-0.3em}{0}}
	\pgfpatharc{-109}{-170}{4\myarrowsize}
	\pgfpatharc{10}{189}{0.58\pgflinewidth and 0.2\pgflinewidth}
	\pgfpatharc{-170}{-115}{4\myarrowsize+\pgflinewidth}
	\pgfpathclose
	\pgfusepathqfillstroke
	\pgfpathmoveto{\pgfpoint{1.5\pgflinewidth-0.3em}{0}}
	\pgfpatharc{109}{170}{4\myarrowsize}
	\pgfpatharc{-10}{-189}{0.58\pgflinewidth and 0.2\pgflinewidth}
	\pgfpatharc{170}{115}{4\myarrowsize+\pgflinewidth}
	\pgfpathclose
	\pgfusepathqfillstroke
	\pgfsetlinewidth{2\pgflinewidth}
}

\pgfarrowsdeclare{cmhook}{cmhook}{
	\pgfsetdash{}{0pt} 
	\pgfsetbeveljoin 
	\pgfsetroundcap 
	\setlength{\myarrowsize}{0.6pt}
	\addtolength{\myarrowsize}{.5\pgflinewidth}
	\pgfarrowsleftextend{-4\myarrowsize-.5\pgflinewidth} 
	\pgfarrowsrightextend{.8\pgflinewidth}
}{
	\setlength{\myarrowsize}{0.6pt} 
  	\addtolength{\myarrowsize}{.5\pgflinewidth}  
 	\pgfsetdash{}{0pt}
	\pgfsetroundcap
	\pgfpathmoveto{\pgfqpoint{0pt}{-4.667\pgflinewidth}}
	\pgfpathcurveto
    {\pgfqpoint{4\pgflinewidth}{-4.667\pgflinewidth}}
    {\pgfqpoint{4\pgflinewidth}{0pt}}
    {\pgfpointorigin}
	\pgfusepathqstroke
}

\newcommand{\ZZ}{\mathbb{Z}}

\DeclareMathOperator{\Jac}{Jac}

\DeclareMathOperator{\Supp}{Supp}

\DeclareMathOperator{\red}{red}

\def\overbar#1#2#3{{%
	\setbox0=\hbox{\displaystyle{#1}}%
	\dimen0=\wd0
	\advance\dimen0 by -#2 
	\vbox {\nointerlineskip \moveright #3 \vbox{\hrule height 0.3pt width \dimen0}%
		\nointerlineskip \vskip 1.5pt \box0}%
}}

\makeatletter
\let\@@seccntformat\@seccntformat
\renewcommand*{\@seccntformat}[1]{%
  \expandafter\ifx\csname @seccntformat@#1\endcsname\relax
    \expandafter\@@seccntformat
  \else
    \expandafter
      \csname @seccntformat@#1\expandafter\endcsname
  \fi
    {#1}%
}
\newcommand*{\@seccntformat@subsection}[1]{%
  \textbf{\csname the#1\endcsname.}
}
\makeatother

\makeatletter
\let\@paragraph\paragraph
\renewcommand*{\paragraph}[1]{%
	\vspace{0.3\baselineskip}%
	\@paragraph{\textbf{#1}}%
}
\makeatother

\counterwithin{equation}{subsection}
\counterwithout{subsection}{section}
\counterwithin{figure}{subsection}

\newtheorem{theorem}[equation]{Theorem}
\newtheorem*{theorem*}{Theorem}
\newtheorem{lemma}[equation]{Lemma}
\newtheorem*{lemma*}{Lemma}
\newtheorem{corollary}[equation]{Corollary}
\newtheorem{proposition}[equation]{Proposition}
\newtheorem*{proposition*}{Proposition}
\newtheorem{conjecture}[equation]{Conjecture}

\theoremstyle{definition}

\newtheorem*{definition*}{Definition}
\theoremstyle{remark}
\newtheorem{remark}[equation]{Remark}

\newtheorem*{example*}{Example}
\newtheorem*{remark*}{Remark}
\newtheorem*{problem*}{Problem}

\newtheorem{notation}{Notation}
\theoremstyle{plain}

\newcommand{\theoremref}[1]{\hyperref[#1]{Theorem~\ref*{#1}}}
\newcommand{\lemmaref}[1]{\hyperref[#1]{Lemma~\ref*{#1}}}
\newcommand{\definitionref}[1]{\hyperref[#1]{Definition~\ref*{#1}}}
\newcommand{\propositionref}[1]{\hyperref[#1]{Proposition~\ref*{#1}}}
\newcommand{\conjectureref}[1]{\hyperref[#1]{Conjecture~\ref*{#1}}}
\newcommand{\corollaryref}[1]{\hyperref[#1]{Corollary~\ref*{#1}}}
\newcommand{\exampleref}[1]{\hyperref[#1]{Example~\ref*{#1}}}

\makeatletter
\let\old@caption\caption
\renewcommand*{\caption}[1]{%
	\setcounter{figure}{\value{equation}}%
	\stepcounter{equation}%
	\old@caption{#1}\relax%
}
\makeatother

\newcounter{intro}

\newtheorem{intro-conjecture}[intro]{Conjecture}
\newtheorem{intro-corollary}[intro]{Corollary}
\newtheorem{intro-theorem}[intro]{Theorem}
\newtheorem{intro-prop}[intro]{Proposition}

\newcommand{\parref}[1]{\hyperref[#1]{\S\ref*{#1}}}

\makeatletter
\newcommand*\if@single[3]{%
  \setbox0\hbox{${\mathaccent"0362{#1}}^H$}%
  \setbox2\hbox{${\mathaccent"0362{\kern0pt#1}}^H$}%
  \ifdim\ht0=\ht2 #3\else #2\fi
  }
\newcommand*\rel@kern[1]{\kern#1\dimexpr\macc@kerna}
\newcommand*\widebar[1]{\@ifnextchar^{{\wide@bar{#1}{0}}}{\wide@bar{#1}{1}}}
\newcommand*\wide@bar[2]{\if@single{#1}{\wide@bar@{#1}{#2}{1}}{\wide@bar@{#1}{#2}{2}}}
\newcommand*\wide@bar@[3]{%
  \begingroup
  \def\mathaccent##1##2{%
    \if#32 \let\macc@nucleus\first@char \fi
    \setbox\z@\hbox{$\macc@style{\macc@nucleus}_{}$}%
    \setbox\tw@\hbox{$\macc@style{\macc@nucleus}{}_{}$}%
    \dimen@\wd\tw@
    \advance\dimen@-\wd\z@
    \divide\dimen@ 3
    \@tempdima\wd\tw@
    \advance\@tempdima-\scriptspace
    \divide\@tempdima 10
    \advance\dimen@-\@tempdima
    \ifdim\dimen@>\z@ \dimen@0pt\fi
    \rel@kern{0.6}\kern-\dimen@
    \if#31
      \overline{\rel@kern{-0.6}\kern\dimen@\macc@nucleus\rel@kern{0.4}\kern\dimen@}%
      \advance\dimen@0.4\dimexpr\macc@kerna
      \let\final@kern#2%
      \ifdim\dimen@<\z@ \let\final@kern1\fi
      \if\final@kern1 \kern-\dimen@\fi
    \else
      \overline{\rel@kern{-0.6}\kern\dimen@#1}%
    \fi
  }%
  \macc@depth\@ne
  \let\math@bgroup\@empty \let\math@egroup\macc@set@skewchar
  \mathsurround\z@ \frozen@everymath{\mathgroup\macc@group\relax}%
  \macc@set@skewchar\relax
  \let\mathaccentV\macc@nested@a
  \if#31
    \macc@nested@a\relax111{#1}%
  \else
    \def\gobble@till@marker##1\endmarker{}%
    \futurelet\first@char\gobble@till@marker#1\endmarker
    \ifcat\noexpand\first@char A\else
      \def\first@char{}%
    \fi
    \macc@nested@a\relax111{\first@char}%
  \fi
  \endgroup
}
\makeatother

\renewcommand{\le}{\leqslant}\renewcommand{\leq}{\leqslant}
\renewcommand{\ge}{\geqslant}\renewcommand{\geq}{\geqslant}

\newcommand{\Q}{\mathbb{Q}}
\newcommand{\OO}{\mathcal{O}}
\newcommand{\be}{\begin{enumerate}}
\newcommand{\ee}{\end{enumerate}}
\newcommand{\bt}{\begin{theorem}}
\newcommand{\et}{\end{theorem}}
\newcommand{\bc}{\begin{corollary}}
\newcommand{\ec}{\end{corollary}}
\newcommand{\bp}{\begin{proof}}
\newcommand{\ep}{\end{proof}}

\newcommand{\ns}{\Sigma}

\newcommand{\Ymu}{{\widetilde{Y}}}
\begin{document}
\keywords{pluricanonical bundles, Fujita's conjecture, effective results.}
\subjclass[2010]{Primary 14C20; Secondary 14F05, 14Q20, 14J17}
\vspace{\baselineskip}

\title[direct images of pluricanonical bundles]{On the effective freeness of the direct images of pluricanonical bundles}

\author[Yajnaseni Dutta]{\bfseries Yajnaseni Dutta}

\address{
Department of Mathematics \\ 
Northwestern University   \\ 
Evanston, IL\\
USA}
\email{ydutta@math.northwestern.edu}

\thanks{}
\setlength{\parskip}{0.19\baselineskip}

\begin{abstract}
We give effective bounds on the generation of pushforwards of log-pluricanonical bundles twisted by ample line bundles. 
This gives a partial answer to a conjecture proposed by Popa and Schnell.  We prove two types of statements: first, 
more in the spirit of the general conjecture, we show generic global generation with predicted bound when the dimesnion 
of the variety if less than 4 and more generally, with a quadratic Angehrn-Siu type bound. 
Secondly, assuming that the relative canonical bundle is relatively semi-ample, we make a very precise statement. 
In particular, when the morphism is smooth, it solves the conjecture with the same bounds, for certain pluricanonical bundles.
\end{abstract}

\maketitle

\subsection{Introduction}  

The main purpose of this paper is to give a partial answer to a version of the Fujita-type conjecture proposed by Popa and Schnell \cite[Conjecture 1.3]{PS}, on the global generation of pushforwards of pluricanonical bundles twisted by ample line bundles. All varieties considered below are over the field of complex numbers.

\noindent
\begin{notation}\label{notation}  We fix \[N=\begin{cases}
n &\text{; when } n\leq 4\\
{n+1 \choose 2}& \text{; otherwise}
\end{cases}\]
in what follows. Our results also work if $N$ was taken to be the 
effective bounds arising from the works of Helmke \cite{Hel97, Hel 99}. 
\end{notation}
\begin{conjecture}[Popa-Schnell]\label{conjPS}
 Let $f : Y \to X$ be a morphism of smooth projective varieties,
with $\dim X=n$, and let $L$ be an ample line bundle on $X$ . Then, for every $k\geq 1$, the sheaf $$f_*\omega_Y^{\otimes k}\otimes L^{\otimes l}$$ is globally generated for $l\geq k(n+1)$.
\end{conjecture}
In \cite{PS}, Popa and Schnell proved the conjecture in the case when $L$ is an ample and globally generated line bundle, and observed that it holds in general when $\dim X = 1$. 
With the additional assumption that $L$ is globally generated, they could use Koll\'ar and Ambro-Fujino type vanishing along with Castelnuovo-Mumford regularity to conclude global generation. We remove the global generation assumption on $L$, by making a 
generation statement at general points with quadratic bounds. 
\begin{intro-theorem}\label{main-b}
Let $f: Y \to X$ be a surjective morphism of projective varieties, with $X$ smooth and $\dim X=n$. Let $L$ be an ample line bundle on $X$. Consider a klt
pair $(Y, \Delta)$ with $\Delta$ a $\Q$-divisor, such that $k(K_Y+\Delta)$ is linearly
equivalent to a Cartier divisor for some integer $k\geq 1$. Denote $P=\OO_Y\big(k(K_Y+\Delta)\big)$. Then the sheaf $$f_*P\otimes L^{\otimes l}$$ is generated by global sections at a general point $x\in X$  for all $l \ge k(N+1)$ with $N$ as in Notation \ref{notation}.

\end{intro-theorem}
\paragraph{A survery of more recent results. }In dimensions higher than $4$, Deng \cite{deng} using analytic extension theorems, showed the
above generation for $l\geq k(n+1)+n^2-n$. More recently Iwai \cite{Iwa17}
 using methods similar to Deng's improved the bound to  $l\ge k(n+1)+\frac{n^2-n}{2}$. Independently at the same time 
the author and Murayama \cite{DM18} improved these results using positivity of $f_*\OO_Y(k(K_Y+\Delta))$, obtaining the bounds similar to Deng, improving the existent results to the log canonical case.

As a particular case of Theorem \ref{main-b}, we have the following corollary, which is a generic version of Conjecture \ref{conjPS} with Angehrn-Siu type bound.
\begin{intro-corollary}
Let $f: Y \to X$ be a surjective morphism of smooth projective varieties, with $\dim X=n$. Let $L$ be an ample line bundle on $X$. Then for all $k \geq 1$, the sheaf $$f_*\omega_Y^{\otimes k} \otimes L^{\otimes l}$$ is generated by global sections at a general point $x\in X$ for all $l$ as in Theorem \ref{main-b}.
\end{intro-corollary}

According to \cite[\S 4]{PS}, this could be interpreted as an effective version of Viehweg's weak-positivity for $f_*\omega_{Y/X}^{\otimes k}$ \cite{Vieh} (also see \cite[Theorem 3.5(i)]{Kol1}).

One can in fact describe the locus on which global generation holds, but not in a very explicit fashion. This suffices however in order to deduce the next Theorem, where assuming semiampleness of the canonical bundle along the smooth fibres, we prove that the global generation holds outside of the 
non-smooth the morphism. 
\begin{intro-theorem}\label{main}
Let $f: Y \to X$ be a surjective morphism of smooth projective varieties, with $\dim X=n$. Suppose $f$ is smooth outside of a 
closed subvariety $B\subset X$. Assume in addition that $\omega_Y^{\otimes k}$ is relatively free outside $B$ for some $k\ge 1$, and let $L$ be an ample line bundle on $X$. Then the sheaf $$f_*\omega_Y^{\otimes k}\otimes L^{\otimes l}$$ is generated by global sections at $x$, for all $x\notin B$ for all $l \geq k(N+1)$.
\end{intro-theorem}
\begin{remark}
Note, for instance, that this applies when $f:Y\to X$ is a projective surjective morphism with generalised Calabi-Yau fibres (i.e. $\omega_F = \OO_F$ for any smooth fibre $F$ of $f$), or with fibres having nef and big canonical bundle (i.e. they are minimal varieties of general type). Indeed, in the second case there is an integer $s\gg 0$ such that $f^*f_*\omega_Y^{\otimes s}\to\omega_Y^{\otimes s}$ is surjective \cite[Theorem 1.3]{Fuj}. 
\end{remark}
In particular, if $f$ is smooth, i.e. $B=\emptyset$, Theorem \ref{main} solves Conjecture \ref{conjPS} for the pluricanonical bundles that are relatively globally generated, however with Angehrn-Siu type bound . 

This in turn leads to an effective vanishing theorem (see Theorem \ref{vanish}), in the case of smooth morphisms, for the pushforwards of pluricanonical bundles that are relatively free. This is in the flavour of \cite[Theorem 1.7]{PS}, but with the global generation assumption on $L$ removed. This vanishing theorem has been improved in \cite{DM18} for $n>4$.

 The proof of Theorem \ref{main-b} is, in part, inspired by arguments in \cite[Theorem 1.4]{PS}. However, since we do not assume that $L$ is globally generated, we need to follow a different path, avoiding Castelnuovo-Mumford regularity. To do this, we need to argue locally around each point and to appeal to the following \textit{local version} of Kawamata's effective freeness result (see \cite[Theorem 1.7]{kaw}), another main source of inspiration for this paper.

\begin{proposition}\label{kaw1}
Let $f: Y \to X$ be a surjective morphism of smooth projective varieties, with $\dim X=n$, such that $f$ is smooth outside of a 
closed subvariety $B$ in $X$. 
Let $\Delta$ be a $\Q$-divisor on $Y$ with simple normal crossing support and coeffecients in $(0,1)$ and let $H$ be a semiample $\Q$-divisor on $X$
such that there is a Cartier divisor $P$ satisfying $$P-(K_Y+\Delta)\sim_{\Q} f^*H.$$ Fixing a point $x\in X\setminus B$,
assume moreover that each strata of $(Y,\Supp(\Delta)$ intersects the fibre above $x$ transversely  or not at all.
Furthermore, let $A$ be a nef and big line bundle on $X$ satisfying $A^n> N^n$ and $A^d.V>N^d$ for any irreducible 
closed subvariety $V \subset X$ of dimension $d$ that contains $x$ and for $N$ as in Notation \ref{notation}. 
Then $$f_*\OO_Y(P)\otimes A$$ is generated by global sections at $x$. 
\end{proposition}
 
\noindent
\begin{remark}\label{kaw-remark}\quad
\be 
\item When $\Delta = 0$ and $B$ is a simple normal crossing divisor, a little more is true: the sheaf $f_*\omega_Y\otimes A$ is in fact generated by global sections 
at every $x \in X$ around which $A$ satisfies Angehrn-Siu type intersection properties. This is Kawamata's freeness result \cite[Theorem 1.7]{kaw}. 
 Kawamata's proof relies on the existence of an effective $\Q$-divisor $D\sim_{\Q} \lambda A$ for some $0<\lambda<1$, such that the pair $(X,D)$ has an isolated log canonical singularity at a given point $x\in X$. Existence of such divisors is known,
 when $A$ satisfies the intersection properties as in the hypothesis of Proposition \ref{kaw1} (see \cite{AS}, \cite[Theorem 5.8]{Ko}). Slightly better bounds are known due to Helmke (\cite{Hel97}, \cite{Hel99}).
Our 
proofs also work with $N$ replaced by Helmke's bounds.

\item The proof proceeds by reducing to the case $\Delta=0$ and then to the situation in Kawamata's result i.e. when $B$ has simple normal crossing support.
We peform the first reduction using an inductive procedure of removing the coeffecients of the components of $\Delta$
via Kawamata coverings \cite[Theorem 4.1.12]{Laz1}. For details see \S \ref{tech-sec}.
\ee
\end{remark}

\noindent
{\bf Acknowledgement.}  I am extremely grateful to my advisor Mihnea Popa for suggesting the problem and for detecting several mistakes in the earlier versions of the proof. I especially thank Lei Wu for several helpful discussions. I would also like to thank Lawrence Ein, Christian Schnell, Valentino Tosatti, Robert Lazarsfeld, Akash Sengupta and Sebasti\'{a}n Olano for several motivating conversations. I
would like to extend my gratitute to
 Takumi Murayama for reading carefully through an earlier version and thoroughly sending me comments.

\subsection{Technical Background}\label{tech-sec}
We begin this section with the proof of the partial generalisation of Kawamata's freeness result. 

\bp[Proof of Proposition \ref{kaw1}]
 Since $H$ is semiample, so is $f^*H$ and therefore by Bertini's theorem we can pick a 
fractional $\Q$ divisor $D\sim_{\Q}f^*H$ with smooth support such that $\Delta+D$ still has simple normal crossings support, 
$\Supp(D)$ is not contained in the support of the $\Delta$ and intersects the fibre over $x$ transversely or not at all and
 $\Delta+D$ still has coeffecitent in $(0,1)$. We rename $\Delta+D$ by $\Delta$. 

Now we proceed by inductively removing the components of $\Delta$. 
\paragraph{Step 1: Kawamata covering of $\Delta$: }  If $\Delta=0$ we move to Step \ref{Step:2}.

Otherwise let $\Delta=\frac{l}{k}D_1+D_2$ with $l,k\in \ZZ_{> 0}$, where $D_1$ is smooth irreducible and $l<k$. 
We choose a Bloch-Gieseker cover $p\colon Z\to Y$ of $Y$ along $D_1$, 
so that $p^*D_1\sim kM$ for some Cartier divisor (possibly non-effective) $M$ on $Z$
and so that the components of $p^*\Delta$ and the fibre $\mathfrak{F}'_x\coloneqq (f\circ p)^{-1}(x)$ 
are smooth and intersect each other transversely or not at all(\cite[Lemma 4.1.11]{Laz1}). 
Moreover since $p$ is flat and $f$ is smooth over a neighbourhood around $x$, 
we can conclude that there is a open neighbourhood $U$ around $x$ such that $f\circ p$ is still smooth over $U$ (\cite[Ex. III.10.2]{Har}.

Set $g=f\circ p$ and denote by $B\subset X$, the branch locus of $g$. Further note that $x\notin B$.

Now, $\omega_Y$ is a direct summand of $p_*\omega_Z$ via the trace map. Therefore $$f_* \OO_Y(P)\otimes A$$ is a direct summand of  
$$g_* \OO_{Z}(K_{Z}+lM+p^*(D_2))\otimes A.$$ 
Hence it is enough to show that the latter	
is generated by global sections at $x$. 

To do this we take the  $k^{\text{th}}$ cyclic cover $q:Y_1\to Z$ of $p^*D_1$. Since $\Supp(p^*\Delta)$ intersects $\mathfrak{F}'_x$ transversely , 
by Lemma \ref{tech} there is an open set $U$ around $x$ such that the components of $p^*\Delta$ intersects all the fibres over $U$ transversely 
or not at all. These intersection properties carries over to $Y_1$ as well.
Further by Lemma \ref{lemma} we see that $g\circ q$ is still smooth over $U$, 
in other words $x$ is not in the branch locus (denoted $B$ again) of $g\circ q$. 
For the ease of notation set  $f_1:=g\circ q$. Note that, 
\[q_*\omega_{Y_1} \simeq \bigoplus_{i=1}^{k-1}\omega_{Z}(p^*D_1-iM) \simeq \bigoplus_{i=1}^{k-1}\omega_{Z}((k-i)M).\]
Indeed, since $p^*D_1 \sim kM$. Further, since $k>l$ the direct sum on the right hand side contains the term $\omega_Z(lM)$ when $i=k-l$.

 Therefore it is enough to show that, 
$$f_{1_*} \OO_{Y_1}(K_{Y_1}+q^*p^*D_2)\otimes A$$ is generated by global sections at $x$.

Proceeding inductively this way, it is enough to show that
$$f_{s_*}\omega_{Y_s}\otimes A$$
is globally generated at $x$, where 
 $f_s:Y_s\to Y$ is the composition of Kawamata covers along the components of $\Delta$ (here $s$ is the number of components of $\Delta$).
We rename $f_s$ by $f$ and $Y_s$ by $Y$. We again call the non-smooth locus of $f_s$ by $B$ and note that $x\notin B$.

\paragraph{Step 2: Base case of the induction: }\label{Step:2}
Take a birational modification $X'$ of $X$ such that $\mu^{-1}(B)_{\red}=:\ns$ in $X'$, as in the diagram below, 
has simple normal crossing support and $X'\setminus \ns \simeq X\setminus B$. In particular, $\mu$ is an isomorphism around $x$.
Let $\tau:Y'\to Y$ be a resolution of the largest irreducible component of the fibre product $X'\times_{X}Y$. We have the following commutative diagram:
\[\begin{tikzcd}
Y' \arrow{r}{\tau} \arrow{d}{h}& Y\arrow{d}{f}\\
X' \arrow{r}{\mu}& X\\
\end{tikzcd}\]
\noindent
Since $\mu$ is an isomorphism over the neighbourhood $U$ around $x$, $\mu^*A$ is a big and nef line bundle
that satisfies the intersection properties, as in the hypothesis, at the point $\mu^{-1}(x)$. 
Moreover since $h$ is smooth outside of the simple normal crossing divisor $\ns$, 
we can apply Kawamata's freeness result \cite[Theorem 1.7]{kaw} to conclude that, 
$h_*\omega_{Y'}\otimes \mu^*A$ is generated by global sections at $\mu^{-1}(x)$. Additionally we have that, 
$$\mu_*h_*\omega_{Y'}\otimes A
\simeq f_*\omega_{Y}\otimes A.$$
Therefore the sheaf $f_*\omega_Y\otimes A$ is generated by global sections at $x$.
	
\ep
\begin{remark}[Local version of Kawamata's theorem]
	When $\Delta=0$, by Szab\'o's Lemma (see e.g. \cite[Theorem 10.45(1)]{Kol13}), we
	can choose $\mu$ in Step \ref{Step:2} of the above proof to be an isomorphism outside the simple
normal crossing locus of $B$ to obtain a local version of Kawamata's theorem. Said differently,
the proof shows that for any morphism $f\colon Y\to X$ between smooth projective varieties, if $x\in X$ has a Zariski neighbourhood $U$ such that $f\colon f^{-1}(U)\to U$ is smooth outside a simple normal crossing divisor then
$f_*\omega_Y\otimes L^{\otimes l}$ is globally generated at $x$, for all $ l\geq N+1$ with $N$ as in Notation \ref{notation}.
\end{remark}

The following two facts were used in the proof of Proposition \ref{kaw1}.

\begin{lemma}\label{tech}
Let $f:Y\to X$ be a smooth and proper surjective morphism of smooth varieties. Let $D$ be a smooth irreducible 
effective divisor on $Y$ such that $D$ intersects a fibre $\mathfrak{F}_x$ over a smooth point $x\in X$ transversely . Then there is an analytic open set $U$ around $x$ such that $f|_{D\cap f^{-1}(U)}:D\cap f^{-1}(U)\to U$ is smooth, or equivalently $D$ intersects all the fibres over $U$ transversely .
\end{lemma}
\bp
For every $y\in D\cap \mathfrak{F}_x$, we choose coordinates $V'_y\subset Y$ around $y$ and coordinates $U_y \subset X$ around $x$ such that $f$ is smooth over $U_y$ and $f(V'_y)\subset U_y$. We choose $V'_y$ such that in these local coordinates we can write $D=(y_m=0)$. Since $D$ intersects the fibre $\mathfrak{F}_x$ over $x$ transversely , we have that
\[\Jac(f|_D)=
  \begin{bmatrix}
\frac{\partial f_1}{\partial y_1} & \ldots &  \frac{\partial f_n}{\partial y_1}\\
 \frac{\partial f_1}{\partial y_2}  & \ldots &\frac{\partial f_n}{\partial y_2}\\
\vdots & \ddots & \vdots\\
\frac{\partial f_1}{\partial y_{m-1}}     &\ldots & \frac{\partial f_n}{\partial y_{m-1}}
\end{bmatrix}_{(y_1,\dots,y_{m-1},0)}
 \]
Let us denote $S_y:= Z( \text{determinant of }n\times n \text{ minors of } f)$. Since $D$ intersects $\mathfrak{F}_x$ transversely , $S_y$ does not intersect $\mathfrak{F}_x$. Then on the neighbourhood $V'_y$, the points where $f$ is smooth is given by the open set $V_y:= V'_y\setminus S_y$ containing $y$. We can cover $\mathfrak{F}_x$ by finitely many such open sets $V_y$. Recall that $f$ is smooth over $\bigcup_y U_y$. 
Pick $$U:=\displaystyle\bigcap_y f(V_y).$$ This is open since a smooth morphism is open and is non-empty since $x\in U$. Then $\Jac(f|_D)$ is non-degenerate on $f^{-1}(U)$.
\ep

Next, we show that the non-smooth locus of a morphism pre-composed with the cyclic cover along a smooth divisor that intersects the smooth fibres of the morphism transversely
remains the same.

\begin{lemma}[Non-smooth locus under cyclic cover]\label{lemma}
Let $f:Y\to X$ be a smooth morphism between smooth varieties. Let $D$ be a divisor on $Y$ that intersects the fibres transversely . Assume in addition that $D\in |L^{\otimes k}|$ for some line bundle $L$ on $Y$ and for some $k\geq 2$. Consider the $k^{th}$-cyclic cover 
$\nu: Y' \to Y$ branched along $D$. Then $f\circ \nu$ is also a smooth morphism.
\end{lemma}
\bp
Pick a point $x\in X$ and a local system of coordinates $x_1, \dots, x_n$ around $x$. Similarly pick a point $y$ on the fibre $\mathfrak{F}_x$ over $x$ and consider a local system of coordinates $y_1, \dots y_m$. In this local coordinates suppose $f$ is given by $(f_1,\dots, f_n)$. Then the Jacobian of $f$ is given by:
\[\Jac(f)=
  \begin{bmatrix}
\frac{\partial f_1}{\partial y_1} & \ldots &  \frac{\partial f_n}{\partial y_1}\\
 \frac{\partial f_1}{\partial y_2}  & \ldots &\frac{\partial f_n}{\partial y_2}\\
\vdots & \ddots & \vdots\\
\frac{\partial f_1}{\partial y_m}     &\ldots & \frac{\partial f_n}{\partial y_m}
\end{bmatrix}
 \]
Since $f$ is smooth, the non-smooth locus of $f$ in $Y$, given by the common zero loci of the determinant of the $n\times n$ minors of $\Jac(f)$, is empty. We denote this locus by $C := Z(\phi_1,\dots,\phi_l)$, where $l={m \choose n}$ and $\phi_i$'s are the determinants of the $n\times n$ minors of $\Jac(f)$.

Assume further that around $y$, $D$ can be written as $D=(y_m=0)$. Then in local coordinates the $k^{th}$-cyclic cover of $D$, $\nu:Y'\to Y$ looks like
\[\nu: (y_1,\dots, y_m) \mapsto (y_1,\dots,y_m^k)\]
Then, 
\[\Jac(f\circ\nu)=
  \begin{bmatrix}
\frac{\partial f_1}{\partial y_1} & \ldots &  \frac{\partial f_n}{\partial y_1}\\
 \frac{\partial f_1}{\partial y_2}  & \ldots &\frac{\partial f_n}{\partial y_2}\\
\vdots & \ddots & \vdots\\
ky_m^{k-1}\frac{\partial f_1}{\partial y_m}     &\ldots & ky_m^{k-1}\frac{\partial f_n}{\partial y_m}
\end{bmatrix}
 \]
As before, the non-smooth locus $C'$ of $f\circ\nu$ in $Y'$ is given by the common zero of the determinant of the $n\times n$ minors of $\Jac(f\circ\nu)$. We write these determinant equations in terms of $\phi_i$'s. 
$$C':= Z(\displaystyle ky_m^{k-1}\phi_1,\dots ,ky_m^{k-1}\phi_s, \phi_{s+1},\dots,\phi_l), $$
 where $\phi_1,\dots,\phi_s$ are equations of the minors that involve the last row of $\Jac(f)$.

We want to show that $C'=\emptyset$. Suppose there is a point $p\in C'$, since $p\notin C$, $p$ must lie on $y_m=0$. Now since $D=(y_m=0)$ intersects all the smooth fibres of $f$ transversely , $p$ cannot be a singular point of $f\big|_D$. In other words, the point $$p\notin Z\big(\phi_{s+1}(y_1,\dots,y_{m-1},0), \dots,\phi_L(y_1,\dots,y_{m-1},0)\big)$$ and hence $p\notin C'$.

\ep

\subsection{Proof of the main theorems}\label{3}
Inspired by \cite{PS}, the strategy is to turn the generation problem for pluricanonical bundles into one for canonical bundles on pairs. We will show that such pair can be carefully chosen to satisfy the properties in the hypothesis of Proposition \ref{kaw1}. 
\bp[Proof of Theorem \ref{main-b}] 

Following the proof of \cite[Theorem 1.7]{PS}, we first take a log resolution $\mu:\Ymu \to Y$ of the base ideal of the adjunction morphism $f^*f_*P \xrightarrow{\pi} P$ and the pair $(Y,\Delta)$. Write:
$$K_{\Ymu} - \mu^*(K_Y+\Delta) = Q - N$$
where	$Q$ and $N$ are effective $\Q$-divisors with simple normal crossing support, with no common components, moreover $N$ has coefficients strictly smaller than 1, and $Q$ is supported on the exceptional locus.
Define:
$$\widetilde{P}:= \mu^*P\otimes \OO_{\Ymu}(\left\lceil Q\right\rceil)$$ and
$$\widetilde{\Delta} := N+\left\lceil Q\right\rceil - Q.$$
Then by definition, the line bundle $\widetilde{P}$ is the same as $ \OO_{\Ymu}\big(k\big(K_{\Ymu}+\widetilde{\Delta}\big)\big)$. Moreover, since $Q$ is exceptional, we have the sheaf isomorphism $\mu_*\widetilde{P} \simeq P$. We rename $\widetilde{Y}$ by $Y$, $\widetilde{P}$ by $P$ and $\widetilde{\Delta}$ by $\Delta$, so that the image of the adjunction morphism $\pi$ is given by $P\otimes \OO_Y(-E)$, for an effective divisor $E$ and so that $Y$ is smooth and the divisor $\Delta+E$ has simple normal crossing support. 

Next, write $\Delta=\sum_ia_i\Delta_i$, where $\Delta_i$'s are the irreducible components of $\Delta$. Let $E_j$'s denote the irreducible components of $E$. Set, $$K=k \times \text{l.c.m. of the denominators of }a_i.$$ Similar to the construction in the proof of Proposition \ref{kaw1}, we take $K^{th}$ Kawamata covers of $\Delta_i$'s and $E_j$'s and denote the composition of these covers by $p:Y'\to Y$. We choose these covers so that $p^*\Delta_i=K\Delta_i'$ and $p^*E_j=KE_j'$ for irreducible divisors $\Delta_i'$ and $E_j'$. We further ensure that $p^*(\Delta+E)$ has simple normal crossing support.

Denote by $B$, the non-smooth locus of $f\circ p$. Consider the following Cartesian diagram: 
\[\begin{tikzcd}
Y\setminus C \arrow{d}{f_V} \simeq V\arrow[hook]{r}{i_V} & Y\arrow{d}{f} \\
X\setminus B =:U \arrow[hook]{r}{i} & X
\end{tikzcd}
\]
where $C=f^{-1}(B)$. 

Fix $x\in X\setminus B$. 
Now, pick a positive integer $m$ which is smallest with the property that the sheaf  
$f_*P\otimes L^{\otimes m}$ is generated by global sections at each point on $U$. Therefore by adjunction,  $P(-E)\otimes f^*L^{\otimes m}$ is globally generated on $V$ and
hence so is $p^*\big(P(-E)\otimes f^*L^{\otimes m}\big)$ on $Y'\setminus \displaystyle p^{-1}(C)$ by the sublinear system $p^*\big|P(-E)\otimes f^*L^{\otimes m}\big|$. By 
Bertini's theorem (see Remark III.10.9.2 \cite{Har} and \cite{Jan}), we can pick $\mathfrak{D}\in \big|P(-E)\otimes f^*L^{\otimes m}\big|$ so that $\mathfrak{D}$ is smooth
outside of $C$ and such that $p^*\mathfrak{D}$ is also smooth outside $p^{-1}(C)$. We further ensure that the divisor $p^*\mathfrak{D}$ intersects the smooth fibre
$\mathfrak{F}'_x$ over $x$ transversely . To simplify notations, we denote $p^{-1}(C)$ by $C$ again. 

We can write:
$$kP+mf^*L\sim \mathfrak{D}+E$$ 

\noindent
From this we get,
$$k\big(K_Y+\Delta\big) \sim_{\Q} K_Y+\Delta+\frac{k-1}{k}\mathfrak{D}+\frac{k-1}{k}E - \frac{k-1}{k}mf^*L,$$
and hence for any integer $l$,
$$k\big(K_Y+\Delta\big)+lf^*L \sim_{\Q} K_Y+\Delta+\frac{k-1}{k}\mathfrak{D}+\frac{k-1}{k}E+\Big(l-\frac{k-1}{k}m\Big)f^*L.$$

Now, since $E$ is the relative base locus of the adjunction morphism $f^*f_*P \xrightarrow{\pi} P$, for every effective Cartier divisor $E'$ such that $E-E'$ is effective we have 
$$f_*\big(P(-E')\big) \simeq f_*P.$$ We would like to pick integral divisors, $E'$ as above so that 
$$\Delta+\frac{k-1}{k}E-E'$$
has coefficients strictly smaller than 1. We do so as follows: write: $$\displaystyle E=\sum_i s_i\Delta_i+\widetilde{E}$$ and $$\displaystyle \Delta = \sum_i a_i\Delta_i$$ such that $\widetilde{E}$ and $\Delta$ do not have any common component. Note that, by hypothesis, $0<a_i< 1$ and $s_i\in \mathbb{Z}_{\ge 0}$. We want to pick non-negative integers $b_i$, 
such that $$0\le a_i+\frac{k-1}{k}s_i - b_i <1$$ and $$b_i\le s_i.$$ Denote by $$\gamma_i:=a_i+\frac{k-1}{k}s_i$$ and note that $\gamma_i< 1+s_i.$
We pick $b_i$ as follows:
 if for some integer $j$, such that $0\le j \le s_i$, we can squeeze $\gamma_i$ between $s_i-j+1> \gamma_i\geq s_i-j$, pick 
$$b_i = s_i-j.$$ Now pick $$E':=\sum_ib_i\Delta_i + \left\lfloor \frac{k-1}{k}\widetilde{E}\right\rfloor.$$ Then we can rewrite the above $\Q$-linear equivalence of divisors as:
$$P-E'+lf^*L \sim_{\Q} K_Y+\widetilde{\Delta}+\frac{k-1}{k}\mathfrak{D}+ \Big(l-\frac{k-1}{k}m\Big)f^*L$$
where $$\widetilde{\Delta}=\Delta+\frac{k-1}{k}E	- E'=\sum_i \alpha_i\widetilde{\Delta}_i.$$ By construction $0<\alpha_i<1$ and $\widetilde{\Delta}$ has simple normal crossing support.

It is now enough to show that the pushforward of the right hand side of the above $\Q$-linear equivalence is globally generated at $x$, for all $l > \displaystyle\frac{k-1}{k}m+N$. Indeed, in that case the left hand side would satisfy similar global generation bounds and by the discussion above
$$f_* P(-E') \otimes L^{\otimes l} \simeq f_* P\otimes L^{\otimes l}.$$

Said differently, this would mean that $$f_*P\otimes L^{\otimes l}$$ is globally generated on $U$ for all $l> \displaystyle\frac{k-1}{k}m+ N$. From our choice of $m$, we must have that $m\le \frac{k-1}{k}m+N+1$. This is the same as $m\le k(N+1)$. As a consequence,$$f_*P\otimes L^{\otimes l}$$ is generated by global sections on $U$ for all $\displaystyle l \geq (k-1)(N+1)+N+1 = k(N+1)$.

It now remains to show that $$f_*\OO_Y\Big(K_Y+\widetilde{\Delta}+\frac{k-1}{k}\mathfrak{D}+ \Big(l-\frac{k-1}{k}m\Big)f^*L\Big)$$ is globally generated at $x$. 
To do so,
we resort to Proposition \ref{kaw1}. 
However the divisor $\widetilde{\Delta}+\frac{k-1}{k}\mathfrak{D}$ may not satisfy the hypothesis of Proposition \ref{kaw1}. 
For instance, it may not have simple normal crossing support. 
Therefore we cannot apply Proposition \ref{kaw1} directly. Since we are only interested in generic global generation though, we can get around these problems. 
The rest of the proof is devoted to this.

We have that $K\alpha_i$ is an integer and by construction, $p$ is a composition of $K^{th}$ Kawamata coverings of the components $\widetilde{\Delta}_i$'s of $\widetilde{\Delta}$. Following the inductive argument as in the proof of Proposition \ref{kaw1}, we see that $$f_*\OO_Y\Big(K_Y+\widetilde{\Delta}+\frac{k-1}{k}\mathfrak{D}+\Big(l-\frac{k-1}{k}m\Big)f^*L\Big)$$ is a direct summand of $$\big(f\circ p\big)_* \OO_{Y'}\Big(K_{Y'}+\frac{k-1}{k}\mathfrak{D}'+(f\circ p)^*\Big(l-\frac{k-1}{k}m\Big)L\Big)$$ where $\mathfrak{D}'=p^*\mathfrak{D}$. Therefore it is enough to show that the latter is globally generated at $x$. 

We are now almost in the situation of Proposition \ref{kaw1}: by our choice of $\mathfrak{D}'$, it intersects the fibre over $x$ transversely , 
however, it may not be klt with simple normal crossing support. 

Note that a log resolution $\mu:Y''\to Y'$ of $\mathfrak{D}'$ is an isomorphism outside $C$. Then write 
$$\mu^*\mathfrak{D}'= D+F$$ where $D$ intersects the fibre over $x$ transversely  and $F$ is supported on $\mu^{-1}(C)$, denoted
by $C$ again. We replace, $Y''$ by $Y'$, rename the divisor $\mu^*\mathfrak{D}'$ by $\mathfrak{D}'$. Therefore, we can assume that $\mathfrak{D}'$ has simple normal crossing support.

To deal with the fact that $F$ may not be klt, consider the effective Cartier divisor $ F'=\left\lfloor \displaystyle\frac{k-1}{k}F \right\rfloor$. 
Since, $\Supp(F')$ is contained in the $C$ and $x\notin B=f(p(C))$, the stalks
\begin{multline*}
\big(f\circ p\big)_* \OO_{Y'}\Big(K_{Y'}+\frac{k-1}{k}\mathfrak{D}'+\Big(l-\frac{k-1}{k}m\Big)(f\circ p)^*L\Big)_x\simeq\\ \big(f\circ p\big)_* \OO_{Y'}\Big(K_{Y'}+\frac{k-1}{k}\mathfrak{D}'-F'+\Big(l-\frac{k-1}{k}m\Big)(f\circ p)^*L\Big)_x
\end{multline*}
are isomorphic. Moreover the global sections of the later embed into the global sections of the former.

Letting 
$$\widetilde{\Delta} :=\displaystyle \frac{k-1}{k}\mathfrak{D}'-F',$$
it is now enough to show that, 
$$\big(f\circ p\big)_* \OO_{Y'}\Big(K_{Y'}+\widetilde{\Delta}+\Big(l-\frac{k-1}{k}m\Big)(f\circ p)^*L\Big) $$
is globally generated at $x$ for $l > \displaystyle\frac{k-1}{k}m+N$.  
The $\Q$-divisor $\widetilde{\Delta}$ satisfies the hypothesis in Proposition \ref{kaw1}. 
Hence the global generation follows from Proposition \ref{kaw1} taking $$H \coloneqq \Big(\left\lceil\frac{k-1}{k}m\right\rceil  - \frac{k-1}{k}m\Big)L$$ 
and $$A\coloneqq \Big(l-\left\lceil\frac{k-1}{k}m\right\rceil\Big)L.$$
Indeed, since $L$ is ample, $H$ is semiample and $A$ satisfies the Angehrn-Siu type intersection properties by the choice of $l$.
\ep 

\begin{remark}
Note that if $\Delta$ already had simple normal crossing support to begin with, then by construction, the loci of generation $U$ in the statement
would contain the largest open set in $X$, over which each strata of $(Y,\Delta+E)$ is smooth.
\end{remark}

The proof of Theorem \ref{main} goes along the same lines. The main difference is that, in this case, we do not start by picking a Kawamata cover, but rather we show directly that, due to the additional relative semi-ampleness assumptions, the above argument works for all $x$ 
outside of the non-smooth locus $B\subset X$ of $f$.

\bp[Proof of Theorem \ref{main}]
As before, we start by replacing $Y$ by a birational modification to assume that the relative base ideal of $\omega_Y^{\otimes k}$ is $\OO_Y(-E)$, for some effective divisor $E$ with simple normal crossing support. Note that, since $\omega_Y^{\otimes k}$ is relatively free over $X\setminus B$, the divisor $E$ is supported on $C:=f^{-1}(B)$. 

Fix a point $x\in X\setminus B$. Consider the following Cartesian diagram:

\[\begin{tikzcd}
Y\setminus C =: V\arrow{d}{f_V} \arrow[hook]{r}{i_V} & Y \arrow{d}{f} \\
X\setminus B=: U\arrow[hook]{r}{i} & X
\end{tikzcd}\]

As in the proof of Theorem \ref{main-b}, we pick a positive integer $m$ which is smallest with the property that the sheaf $f_*\omega_Y^{\otimes k}\otimes L^{\otimes m}$ is generated by global sections at each point of $U$. Then $f^*f_*\big(\omega_Y^{\otimes k}\otimes f^*L^{\otimes m}\big)$ is also generated by global sections on $V$. Therefore by adjunction, so is $\omega_Y^{\otimes k}(-E)\otimes f^*L^{\otimes m}$ on $V$. As a consequence, we can pick a divisor $\mathfrak{D}\in \big|\omega_Y^{\otimes k}(-E)\otimes f^*L^{\otimes m}\big|$ such that $\mathfrak{D}$ is smooth outside of $C$ and intersects the fibre $\mathfrak{F}_x$ over $x$ transversely . 

After replacing $Y$ with a birational modification that is an isomorphism outside of $C$, we may assume that $\mathfrak{D}=D+F$, where $D$ is smooth, intersects the fibre $\mathfrak{F}_x$ over $x$ transversely  and does not share any component with $E$. Moreover, we assume that the divisor $F$ is supported on $C$, and that $E+D+F$ has simple normal crossing support. 

Write
\[kK_{Y}+mf^*L \sim D+F+E.\]
From this we can write,
\[kK_Y+lf^*L\sim_{\Q}K_Y+\frac{k-1}{k}D+\frac{k-1}{k}(F+E)+\Big(l-\frac{k-1}{k}m\Big)f^*L\]
for any integer $l$.
 
Now consider the effective divisor $ \left\lfloor \frac{k-1}{k}(E+F)\right\rfloor$ and denote the fractional part by $$D':=\displaystyle\frac{k-1}{k}(E+F)-\left\lfloor \frac{k-1}{k}(E+F)\right\rfloor.$$ We obtain the following $\Q$-linear equivalence:
\[kK_Y-\left\lfloor \frac{k-1}{k}(E+F)\right\rfloor+lf^*L\sim_{\Q}K_Y+\frac{k-1}{k}D+D'+\Big(l-\frac{k-1}{k}m\Big)f^*L.\]
Denote $\Delta:=\displaystyle\frac{k-1}{k}D+D'$. It is now enough to show that $$f_*\OO_Y\Big(K_Y+\Delta+\big(l-\frac{k-1}{k}m\big)f^*L\Big)$$ is generated by global sections at $x$ for all $l> \frac{k-1}{k}m+N$. Indeed, this would imply that the left hand side of the equation also satisfies similar global generation bounds, i.e. $f_*\OO_Y\Big(kK_Y- \left\lfloor \frac{k-1}{k}(E+F)\right\rfloor+lf^*L\Big)$ is globally generated at $x$ for all $ l> \frac{k-1}{k}m+ N$.  But note that the divisor $E+F$ is supported on $C$ and $x\notin B$. Therefore the stalks 
$$f_*\OO_Y\Big(kK_Y-\left\lfloor \frac{k-1}{k}(E+F)\right\rfloor+lf^*L\Big)_x \simeq f_*\OO_Y\big(kK_Y+lf^*L\big)_x$$ 
are isomorphic. Moreover the global sections of the former embeds into the global sections of the later. Said differently, this would in turn imply that $$f_*\OO_Y\big(kK_Y+lf^*L\big)$$ is globally generated on $U$ for all $l> \displaystyle\frac{k-1}{k}m+ N$. But from our choice of $m$ and from similar arguments as in the proof of Theorem \ref{main-b}, it follows that for all $\displaystyle l \geq  k(N+1),$ the sheaf $$f_*\omega_Y^{\otimes k}\otimes L^{\otimes l}$$ is generated by global sections on $U$.\\

It now remains to show that $$f_*\OO_Y\Big(K_Y+\Delta+\big(l-\frac{k-1}{k}m\big)f^*L\Big)$$ is generated by global sections at $x$ when $l-\frac{k-1}{k}m> N$. 
But this follows from Proposition \ref{kaw1} taking $$H \coloneqq (\left\lceil\frac{k-1}{k}m\right\rceil  - \frac{k-1}{k}m)L$$ 
and $$A\coloneqq \Big(l-\left\lceil\frac{k-1}{k}m\right\rceil\Big)L.$$
Indeed, since $L$ is ample, $H$ is semiample and $A$ satisfies the Angehrn-Siu type intersection properties by the choice of $l$.
Moreover the divisor $\Delta = \frac{k-1}{k}D+D'$ is klt and has simple normal crossing support, its components was chosen to intersect the fibre 
$\mathfrak{F}_x$ transversely . 
\ep

\subsection{An effective vanishing theorem}
When $f$ is smooth, Koll\'ar's vanishing theorem applied to the right hand side of the equivalence, leads to the following vanishing statement for pluricanonical bundles, with essentially the same proof.
\bt[Effective Vanishing Theorem]\label{vanish}
Let $f: Y \to X$ be a smooth surjective morphism of smooth projective varieties, with $\dim X=n$. Assume in addition that $\omega_Y^{\otimes k}$ is relatively free for some $k\geq 1$, and let $L$ be an ample line bundle on $X$. Then, 
$$H^i\Big(X, f_*\omega_Y^{\otimes k}\otimes L^{\otimes l}\Big) = 0$$
for all $i>0$ and $l\geq\displaystyle k(N+1)-N$ with $N$ as in Notation \ref{notation}.
\et
\bp
Since $f$ is smooth, by Theorem \ref{main}, we know that the sheaf $f_*\omega_Y^{\otimes k}\otimes L^{\otimes l}$ is globally generated for all $l\geq k(N+1)$. Therefore by adjunction $\omega_Y^{\otimes k}\otimes f^*L^{\otimes k(N+1)}$ is globally generated as well. As a consequence, we can pick a smooth divisor $D\in \big|\omega_Y^{\otimes k}\otimes f^*L^{\otimes k(N+1)}\big|$ such that $D$ intersects the fibre $\mathfrak{F}_x$ over $x$ transversely . 

Write:
$$kK_Y+k(N+1)f^*L \sim D.$$
This is the same as 
$$kK_Y+lf^*L \sim_{\Q} K_Y+ \frac{k-1}{k}D+\big(l-(k-1)(N+1)\big)f^*L,$$
for any integer $l$.
By applying Koll\'ar's vanishing theorem \cite[Corollary 10.15]{Kol2} on the right hand side, we get that 
$$H^i\Big(X, f_*\OO_Y\big(K_Y+ \frac{k-1}{k}D+\big(l-(k-1)(N+1)\big)f^*L\big)\Big) = 0$$
for all $i>0$ and $l>\displaystyle (k-1)(N+1)$. Therefore, the left hand side satisfies similar vanishing properties
$$H^i\Big(X, f_*\omega_Y^{\otimes k}\otimes L^{\otimes l}\Big) = 0$$
for all $i>0$ and $l\geq\displaystyle k(N+1)-N$
\ep

\begin{remark}
The above bound is replaced in \cite[Theorem 5.3]{DM18} by $k(n+1)-n$ for all $n$. This is an improvement for $n\geq 4$.
\end{remark}

\end{document}